\documentclass{article}
\usepackage[usenames,dvipsnames]{pstricks}
\usepackage{epsfig}
\usepackage{pst-grad} 
\usepackage{pst-plot} 
\usepackage{amsmath,amsthm,amssymb}
\usepackage{pst-plot,color}
\usepackage{pst-grad}
\usepackage{epsfig}
\usepackage[usenames,dvipsnames]{pstricks}
\newcommand{\ds}{\displaystyle}
\newcommand{\lp}{\mbox{L}}
\begin{document}
\large
\title{Convergence and divergence of averages along subsequences in certain Orlicz Spaces}
\author{Christopher Wedrychowicz \\ SUNY Albany\\ Albany, NY, 12222, USA}
\maketitle
\begin{abstract}
The classical theorem of Birkhoff states that the $\ds T^{N}f(x)=\frac{1}{N}\sum_{k=0}^{N-1}f\left(\sigma^{k}x\right)$ converges almost everywhere 
for $\ds x\in X$ and $\ds f\in L^{1}(X)$, where $\ds \sigma$ is a measure preserving transformation of a probability measure space $\ds X$. It was shown that there are operators 
of the form $\ds T^{N}f(x)=\frac{1}{N}\sum_{k=0}^{N-1}f\left(\sigma^{n_{k}}x\right)$ for a subsequence $\ds \{n_{k}\}$ of the positive integers that converge in some $\ds L^{p}$ spaces 
while diverging in others. The topic of this talk will examine this phenomenon in the class of Orlicz spaces $\ds \left\{L\mbox{Log}^{\beta}L:\beta>0\right\}$.
\end{abstract}
\section{}
\theoremstyle{definition}
\newtheorem{dynsys}{Definition}[section]
\begin{dynsys}
Let $\ds (X,\mathcal{B},\mu)$ be a measure space. Let $\ds T:X\rightarrow X$ be a 
one-to-one, onto map such that $\ds \mu(T^{-1}A)=\mu(A)$ $\ds \forall A\in \mathcal{B}$. Then $T$ is called a measure preserving transformation and $\ds (X,\mathcal{B},\mu,T)$ is called a dynamical system.
\end{dynsys}
\theoremstyle{remark}
\newtheorem{eg}{Example}[section]
\begin{eg}An example of central importance to this work is when $\ds X=[0,1)$, $\ds \mu$ is Lebesgue measure, and $\ds \mathcal{B}$ is the $\ds \sigma-$algebra of Borel sets and $T$ is defined by $\ds T(x)= x+ \alpha\mbox{mod}(1)$ where $\ds\alpha \in [0,1)$. It is equivalent to realizing $\ds [0,1)$ as the unit circle and T as a rotation by $\ds 2\pi \alpha$.
\end{eg}

\theoremstyle{plain}
\newtheorem{changevar}{Theorem}[section]
\begin{changevar}[Change of Variable Formula]
Let $\ds (X,\mathcal{B},\mu)$ , $\ds (Y,\mathcal{C},\lambda)$ be measure spaces and let $\ds \Phi: X \rightarrow Y$ be a measurable map in the sense that $\ds \Phi^{-1}(A) \in \mathcal{B}$ for all $\ds A \in \mathcal{C}$. Then ,
\begin{eqnarray*}
\int_{X}f(\Phi(x))\mbox{d}x = \int_{Y} f(y)\mbox{d}y .
\end{eqnarray*} 
\end{changevar}
\theoremstyle{definition}
\newtheorem{erg}[dynsys]{Definition}
\begin{erg}
If the averages 
\begin{eqnarray*}
\frac{1}{N}\sum_{k=0}^{N-1}\mu(T^{-k}A\cap B)\xrightarrow{N\rightarrow\infty}\mu(A)\mu(B) & \quad \forall A,B\in\mathcal{B} ,
\end{eqnarray*}
Then $T$ is called ergodic.
\end{erg}
\theoremstyle{remark}
\newtheorem{egerg}[eg]{Example}
\begin{egerg}
If $\ds\alpha$ is irrational then $T$ is ergodic as defined in the previous example; if $\ds\alpha$ is rational then $T$ is not ergodic.
\end{egerg}
A theorem of fundemental importance in ergodic theory is Birkhoff's Theorem, which is stated as follows,
\theoremstyle{plain}
\newtheorem{birk}[changevar]{Theorem}
\begin{birk}[Birkhoff]
Let $\ds (X,\mathcal{B},\mu,T)$ be a dynamical system and $\ds (X,\mathcal{B},\mu)$ be a $\sigma-$finite measure space then,
\begin{eqnarray*}
\frac{1}{N}\sum_{k=0}^{N-1}f(T^{k}x)\xrightarrow{N\rightarrow\infty}E(f |J)(x) \,\mbox{a.e.}
\end{eqnarray*}
where $J$ is the $\sigma-$algebra of invariant sets.
\end{birk}
There have been many attempts to generalize Birkhoff's \\
Theorem. One in particular is connected to the topic of this thesis.  \\
Let $\ds \{n_{k}\}$ be an increasing sequence of positive integers. One may ask the following question:\\
Do the averages,
\begin{eqnarray*}
\frac{1}{N}\sum_{k=0}^{N-1}f(T^{n_k}x) 
\end{eqnarray*}
converge a.e $\ds \forall f$ in some subspace of $\lp_{1}$?\\
Much work has been done in this area. For example when $\ds n_{k}=k^{2}$ Bourgain has shown that the averages converge a.e. $\ds \forall f\in \lp_{p}$ where $\ds p>1$. The problem of the case $\ds p=1$ remained open for some time. Recently it was shown that for every dynamical system there exists a function $\ds f\in \lp_{1}$ such that the averages do not converge a.e.  The following question dealing with subsequences leads to the topic of this thesis. \\ \vspace{.0001in} \\
First a few definitions,
\theoremstyle{definition}
\newtheorem{unigood}[dynsys]{Definition}
\begin{unigood}
An increasing sequence of integers $\ds (n_{k})$ is called universally $\ds \lp_{p}$ good,
if the averages
\begin{eqnarray*}
 A_{N}(x)=\frac{1}{N}\sum_{k=0}^{N-1}f(T^{n_{k}}x) 
\end{eqnarray*} 
converge a.e. for $\ds x\in X$, $\ds\forall f\in \lp_{p}$ and for all dynamical systems $\ds (X,\mathcal{B},\mu,T)$.\\
 A sequence is called $\ds \lp_{p}$ universally bad if for every dynamical system there exists a function $\ds f\in \lp_{p}$ such that the averages $\ds A_{N}$ fail to converge a.e.
\end{unigood}
\begin{itemize}  
\item[Question 1:]
Does there exist an increasing sequence of integers $\ds (n_{k})$ that is $\ds \lp_{p}$ universally good while $\ds \lp_{q}$ universally bad for all $\ds q<p$?  
\item[Question 2:]
Does there exist a sequence that is $\ds L_{p}$ universally bad but $\ds \lp_{q}$ universally good for all $\ds q\geq p$?
\end{itemize}
The first question was answered affirmatively by Reinhold while the second was answered affirmatively by Bellow.
\theoremstyle{definition}
\newtheorem{llogl}[dynsys]{Definition}
\begin{llogl}
The space of functions $\ds \lp^{s}\mbox{Log}^{p}\lp$ is defined as 
\begin{eqnarray*}
\lp^{s}\mbox{Log}^{p}\lp=\{f\in \lp_{1}: \int \left|f^{s}(x)\right|\mbox{Log}^{p}(\left|f(x)\right|+1)\mbox{d}x< \infty\}
\end{eqnarray*}
\end{llogl}
The notions of universally good and bad extend to the above spaces in an obvious way.
\theoremstyle{remark}
\newtheorem*{q}{Question}
\begin{q}
Given $\ds p$ and a dynamical system $(X,\mathcal{B},\mu,T)$ does there exist an increasing sequence of integers $\ds (n_{k})$ such that the averages $\ds A_{N}(x)$ converge a.e. for all $\ds f\in \lp\mbox{Log}^{q}\lp$ with $\ds q>p$ while there exists a function $f\in \lp\mbox{Log}^{p}\lp$ such that the averages $\ds A_{N}(x)$ fail to converge a.e.
\end{q}
When $p>1$ the answer will be affirmative while for $p\leq 1$ we will prove that there exists an increasing sequence $(n_{k})$ such that $A_{N}(x)$ converge a.e for all $f\in \lp\mbox{Log}^{q}\lp$ for all $\ds q>p$, for $\ds q<p$ there exists a function $f\in \lp\mbox{Log}^{q}\lp$ such that $A_{N}$ fail to converge a.e. while the behavior of $A_{N}(x)$ for functions in $\lp\mbox{Log}^{p}\lp$ is unknown.\\
\section{}
\theoremstyle{plain}
\newtheorem{banach}{Theorem}[section]
\begin{banach}[Banach's principle]
If $T^{\ast}f(x) <\infty$ a.e. for all $f\in B$ where B is a Banach space of functions contained in $\ds\mbox{L}_{1}$ then there is a positive, decreasing function $\ds C(\lambda)$ defined for $\ds \lambda>0$ that goes to zero as $\ds\lambda\rightarrow \infty$ such that for all $\ds f\in B$ we have 
\begin{eqnarray*}
\mu\{x:T^{\ast}f(x)>\lambda\left\|f\right\|_{\mathcal{B}}\} \leq C(\lambda) .
\end{eqnarray*} 
\end{banach}
\newtheorem{ban}[banach]{Theorem}
\begin{ban}
Let $\ds (X,\mathcal{B},\mu)$ be a probability space and $\ds S\subseteq \lp_{1}$ be a Banach space. 
If $\ds \{T_{n}\}$ is a sequence of bounded operators such that 
\begin{eqnarray*}
T^{\ast}f(x)=\sup_{n}T_{n}f(x)<\infty \quad \mbox{a.e.}
\end{eqnarray*}
for every $\ds f\in S$ then the set of functions in $S$ such that $\ds T_{n}f(x)$ converges a.e. is closed.
\end{ban}
In order to establish the inequality above one often establishes a weak \\
maximal inequality for the sublinear operator $\ds T^{\ast}$, that is an inequality of the form
\begin{eqnarray*}
\mu(\{x:T^{\ast}(x)\geq \lambda \})\leq C(\lambda)
\end{eqnarray*}
where $\ds C(\lambda)$ is a monotone decreasing function such that \\
\begin{eqnarray*}
C(\lambda)\xrightarrow{\lambda\rightarrow\infty}0 .
\end{eqnarray*}
\theoremstyle{definition}
\newtheorem{con}[dynsys]{Definition}
\begin{con}
Let $\ds\Phi(x)$ be a function such that
\begin{enumerate}
\item $\ds \Phi$ is continuous and convex.
\item $\ds \Phi(x)=\Phi(-x)$ 
\item 
$\ds \frac{\Phi(x)}{x} \xrightarrow{x\rightarrow 0}   0 $\\ \vspace{.01in} 
\item
$\ds  \frac{\Phi(x)}{x} \xrightarrow{x\rightarrow\infty} \infty$
\end{enumerate}
Let $\ds\lp_{\Phi}=\left\{f\in \lp_{1}: \int \Phi(f(x))\mbox{d}x<\infty\right\}$ \\
then $\ds\lp_{\Phi}$ is a Banach space under the following norm
\begin{eqnarray*}
\left\|f\right\|_{\Phi}=\inf\left\{k:\int \Phi\left(\frac{f}{k}\right)\mbox{d}x <1\right\}
\end{eqnarray*}
\end{con}
\theoremstyle{plain}
\newtheorem{thmdegu}[banach]{Theorem}
\begin{thmdegu}[Sawyer's Theorem]
Let $\ds (X,\mathcal{B},\mu)$ be a probability measure space.  Let $\ds \{T_{k}\}$ be a sequence of positive linear operators from $\ds\mbox{L}_{\Phi}$ to the set of measureable functions on $\ds X$.  Assume that the $\ds\{T_{k}\}$'s commute with a family $\{S\alpha\}$ of measure preserving maps from $X$ to $X$ that mix the measurable sets of $X$.  Assume further that the function $\ds\Phi$ satisfies the following:
\begin{eqnarray*}
\mbox{If}\quad y\geq 1\, ,\, x\geq \frac{1}{y}\quad \mbox{then}\quad \Phi(xy)\leq C(\Phi(y))^{p}\Phi(x) .
\end{eqnarray*}  
Then the following are equivalent:
\begin{enumerate}
\item $T^{\ast}$ satisfies an inequality of the form 
\begin{eqnarray*}
\mu\{x:T^{\ast}f(x)\geq\lambda\} \leq C\int \Phi\left(\frac{f}{\lambda}\right) .
\end{eqnarray*}
\item For each $\ds f \in \mbox{L}_{\Phi}$ , $T^{\ast}f(x) < \infty$ .
\end{enumerate}
\end{thmdegu}
\begin{proof} The following lemma is of central importance to the proof of the 
theorem.
\newtheorem{sawyerlm}[banach]{Lemma}
\begin{sawyerlm}
Let $(X,B,\mu)$ be a probability measure space.  Let \\
$S_{\alpha}:X\rightarrow X$ be a collection of measure preserving maps that mix the measurable sets of $X$. Then if $\{A_{k}\}$ is a sequence of measurable sets of $X$ such that $\sum \mu(A_{k}) = \infty$, there exists a sequence \\ $\{S_{k}\} \subseteq (S_{\alpha})$  such that almost every $x\in X$ is in infinitely many of the sets $S_{k}^{-1}(A_{k})$.
\end{sawyerlm}
Assume that $T^{\ast}$ does not satisfy an inequality of the form
\begin{eqnarray*}
\mu(\{x: T^{\ast}f(x) \geq \lambda \}) \leq C\int \Phi\left(\frac{f}{\lambda}\right) .
\end{eqnarray*}  
Then fix a sequence $c_{k}$ increasing to infinity, $c_{k}>0$. Then there exists a sequence \\
$\ds\{f_{k}\}\subseteq L_{\Phi}$, $\lambda_{k}>0$ such that,
\begin{eqnarray*}
\mu\{T^{\ast}f(x)\geq \lambda_{k}\}> c_{k}\int \Phi\left(\frac{f_{k}}{\lambda_{k}}\right) .
\end{eqnarray*}  
Call $\ds g_{k}=\frac{f_{k}}{\lambda_{k}}$, $\ds A_{k}=\{g_{k}\geq 1\}$. Then,\\

\begin{eqnarray*}
1\geq \mu(A_{k}) \geq c_{k} \int\Phi(g_{k}) .
\end{eqnarray*}
Let $h_{k}$ be natural numbers such that $\ds 1\leq h_{k} \mu(A_{k}) \leq 2$ and take $h_{k}$ copies of $A_{k}$ denoted by $\ds A_{k}^{1},\cdots,A_{k}^{h_{k}}$.  Thus $\ds\sum_{k=1}^{\infty} \sum_{j=1}^ {h_{k}} \mu(A_{k}^{j}) = \infty$ \\
and by the previous lemma there are $\ds S_{k}^{j} \in (S_{\alpha})$ such that almost every $x\in X$ is in infinitely many of the sets $\ds(S_{k}^{j})^{-1}(A_{k}^{j}).$\\
Define a function 
\begin{eqnarray*}
F(x)= \sup_{\ds k\geq1\atop \ds 1\leq j \leq h_{k}} \alpha_{k}S_{k}^{j}g_{k}^{j}(x)
\end{eqnarray*}
where $\ds g_{k}^{j}=g_{k}$ and the constants $\alpha_{k}$ will be determined later.\\
We have 
\begin{eqnarray*}
F(x)=P(x) + Q(x)
\end{eqnarray*}
where 
\begin{eqnarray*}
P(x)= \sup_{\ds k\geq1\atop\ds 1\leq j\leq h_{k}} S_{k}^{j}g_{k}^{j}\geq \frac{1}{\alpha_{k}}
\end{eqnarray*} 
and $\ds Q(x)$ is a function bounded by $1$ .\\
Then 
\begin{eqnarray*}
 P(x) & \leq & \sum_{k,j \in R}\Phi(\alpha_{k}S_{k}^{j}g_{k}^{j}(x)) \\
      & \leq &\sum_{k,j\in R} C[\Phi(\alpha_{k})]^{p}\Phi(S_{k}^{j}g_{k}^{j}(x))\\
      & \leq & C\sum_{k=1}^{\infty}[\Phi(\alpha_{k})]^{p} \sum_{j=1}^{h_{k}}\Phi(S_{k}^{j}g_{k}^{j}(x)) .
 \end{eqnarray*}
and so 
\begin{eqnarray*}
\int \Phi(P(x)) & \leq & C\sum_{k=1}^{\infty} [\Phi(\alpha_{k}]^{p}h_{k} \int \Phi(g_{k}) \\
                & \leq & C\sum_{k=1}^{\infty} [\Phi(\alpha_{k})]^{p}\frac{ \mu(A_{k})}{c_{k}}h_{k} \\
                & \leq & C\sum_{k=1}^{\infty}\frac {[\Phi(\alpha_{k})]^{p}}{c_{k}} ,
\end{eqnarray*} 
by the change of variable formula for measure preserving transformations.
Given that the sequence $\ds\left\{\frac{1}{c_{k}}\right\}$ sums, the $\ds\{\alpha_{k}\}$ may be chosen so that the above sum is finite and the $\alpha_{k}$ increase to infinity.  The remainder of the argument is the same as in ~\cite{degu}.
\end{proof}
If we fix a sequence $\ds (n_{k})$ to each dynamical system, we may \\
associate a constant $\ds C(n_{k})$ such that 
\begin{eqnarray*}
 \mu\{x:T^{\ast} > \lambda \} \leq C(\{n_{k}\}) \int \Phi\left(\frac{f}{\lambda}\right) .
\end{eqnarray*}
We may then consider the minimal such constant so that a similar inequality holds in all dynamical systems.  The so-called Conze's principle asserts a condition in which we may conclude that such a minimal constant exists and is finite. As a result this will by Sawyer's Theorem confirm whether a sequence is universally good or not.
\theoremstyle{plain}
\newtheorem{conze}[banach]{Theorem}
\begin{conze}[Conze's Principle]
For a given sequence $\ds (n_{k})$ to have it's associated minimal constant finite, it is enough that there exists a single ergodic dynamical system $\ds(X,\mathcal{B},\mu,T)$ such that the averages
\begin{eqnarray*}
\frac{1}{N}\sum_{k=0}^{N-1} f(T^{n_{k}}) \quad\mbox{converge}\quad\mbox{a.e.}
\end{eqnarray*}
\end{conze}
\section{}
The main candidates for such sequences will be perturbations of block sequences.
A block of integers is a set of the form $\ds B=[n,n+1,\cdots,n+k-1]$ of consecutive integers. We will let $\ds |B|=k$ denote the number of integers in $B$ and will refer to it as the length of $B$. A block sequence is a sequence $\ds \{n_{k}\}$ that can be arranged into blocks $\ds B_{1},B_{2},\cdots$ as a set $\ds \{n_{k}\}=\bigcup_{k=1}^{\infty}B_{k}$.
Let $\ds D_{k}$ be an arbitrary collection of integers between $\ds B_{k}$ and $\ds B_{k+1}$.
The collection $\ds \bigcup_{k=1}^{\infty}D_{k}$ will be referred to as a perturbation of the block sequence $\ds \bigcup_{k=1}^{\infty}B_{k}$ and the  sequence whose elements are $\ds \bigcup_{k=1}^{\infty}B_{k}\cup D_{k}$ will be referred to as a perturbed block sequence.
The following theorem is a generalization of a Theorem of Bellow. It essentially states that if we begin with a block sequence, which is uiniversally good in a certain subspace of $\ds \lp_{1}$ there is a certain degree to which we may perturb it so that the resulting sequence is also universally good in that subspace. 
\theoremstyle{plain}
\newtheorem{reinhold}{Theorem}[section]
\begin{reinhold}[Reinhold]
Let $\ds B_k$ and $\ds D_k$ be a block sequence and a perturbation of that block sequence. If the sequence $\ds\bigcup_{k=1}^{\infty}B_k$ is universally good for $\ds\mbox{L}_{\infty}$ and
\begin{eqnarray*}
\frac{d_1+\cdots+d_k}{l_1+\cdots+l_k}<\infty
\end{eqnarray*}
then the sequence $\ds\bigcup_{k=1}^{\infty}B_k\cup D_k$ is also universally good for $\ds\mbox{L}_{\infty}$ .
\end{reinhold}
\newtheorem{first}[reinhold]{Theorem}
\begin{first}
Let $\ds \bigcup_{k=1}^{\infty} B_{k}$ be a block sequence that is universally good in the Orlicz space $\ds\lp_{\Phi}$, and let $\ds \bigcup D_{k}$ be as above. Then if
\begin{eqnarray*}
\sum_{k=1}^{\infty}\frac{1}{\Phi\left(\frac{l_{1}+\cdots+l_{k}}{d_{1}+\cdots+d_{k}}\right)}<\infty
\end{eqnarray*}
then the sequence $\ds \bigcup_{k=1}^{\infty} B_{k}\cup D_{k}$ is also universally good in $\ds\lp_{\Phi}$
\end{first}
\begin{proof}
We proceed as in ~\cite{krl}. \\
Let 
\begin{eqnarray*}
C & = & \bigcup B_{k}\cup D_{k} ,\\
b_{n} & = & \left| \bigcup_{k=1}^{\infty} B_{k} \cap [0,n]\right|\:\mbox{and} \\ 
c_{n} & = & \left| \bigcup_{k=1}^{\infty} D_{k} \cap [0,n]\right| . 
\end{eqnarray*}
The averages $\ds A_{n}f(x)=\frac{1}{\left|C\cap[0,n]\right|}\sum_{u\in C\cap [0,n]}f(T^{u}x)$\\
can be written as the convex combination\\
\begin{eqnarray*}
A_{n}f(x) & = & \frac{b_{n}}{b_{n}+c_{n}}\left(\frac{1}{b_{n}} \sum_{u\in \bigcup_{k=1}^{\infty} {B_{k} \cap [0,n]}} f(T^{u}x)\right)+\frac{c_{n}}{b_{n}+c_{n}}\left(\frac{1}{c_{n}}\sum_{u\in \bigcup_{k=1}^{\infty} D_{k}\cap[0,n]}f(T^{u}x)\right) \\
& = & \frac{b_{n}}{b_{n}+c_{n}}A_{n}^{B}f(x) + \frac{c_{n}}{b_{n}+c_{n}}A_{n}^{D}f(x).
\end{eqnarray*}  
To establish a.e. convergence it is enough to do so on each piece sepparately.\\
First we observe that since \\
\begin{eqnarray*}
\frac{1}{\Phi\left(\frac{l_{1}+\cdots+l_{k}}{d_{1}+\cdots+d_{k}}\right)} \rightarrow 0 ,
\end{eqnarray*}
we have 
\begin{eqnarray*}
\Phi\left(\frac{l_{1}+\cdots+l_{k}}{d_{1}+\cdots+d_{k}}\right)\rightarrow \infty 
\end{eqnarray*}
so 
\begin{eqnarray*}
\frac{l_{1}+\cdots+l_{k}}{d_{1}+\cdots+d_{k}} \rightarrow \infty 
\end{eqnarray*}
and hence it's reciprocal goes to $0$.\\
This implies by the previously stated theorem that the averages of functions in $\ds\lp_{\infty}$ converge a.e. \\
We have
\begin{eqnarray*}
 \frac{c_{n}}{b_{n}}=\left\{
\begin{array}{ll}
\ds\frac{d_1+\cdots+d_{k-1}}{l_1+\cdots+l_{k-1}+s_k} &  \mbox{if}\: k \: \mbox{is}\:\mbox{the}\:\mbox{smallest}\:\mbox{integer}\\ 
                                                      &  \mbox{such}\:\mbox{that}\:  B_k \:\mbox{is}\:\mbox{not}\:\mbox{contained} \\
                                                      &  \mbox{in}\: [0,n] \:\mbox{and}\: B_k\cap[0,n]\neq\emptyset , \\
\ds\frac{d_{1}+ \cdots +d_{k-2}+r_{k-1}}{l_{1}+\cdots +l_{k-1}} & \mbox{if}\: k \:\mbox{is}\:\mbox{the}\:\mbox{smallest}\:\mbox{integer}\\
                                                                &  \mbox{such}\:\mbox{that}\: B_k \:\mbox{is}\:\mbox{not}\:\mbox{contained} \\
                                                                &  \mbox{in}\: [0,n] ,\: B_k\cap[0,n]=\emptyset \:\mbox{and}\: B_{k-1}\subset[0,n]
\end{array}\right.\end{eqnarray*}  \\ \vspace{.01in} \\
where $\ds 0\leq r_{k-1}\leq d_{k-1}$ and $\ds 0\leq s_{k} \leq l_{k}$. \\ \vspace{.01in} \\
In either case 
\begin{eqnarray*}
 \frac{c_{n}}{b_{n}} \leq\frac{d_{1}+\cdots+d_{k-1}}{l_{1}+\cdots+l_{k-1}}\rightarrow 0 ,
\end{eqnarray*}
and so
\begin{eqnarray*}
\frac{b_{n}}{b_{n}+c_{n}} \rightarrow 1 . 
\end{eqnarray*}
Therefore
\begin{eqnarray*}
\frac{b_{n}}{b_{n}+c_{n}} A_{n}^{B}f(x)
\end{eqnarray*}
 converges a.e. since\\ $\ds\bigcup_{k=1}^{\infty}B_{k}$ is universally good in $\ds\lp_{\Phi}$.\\
Consider the following operator:\\
\begin{eqnarray*}
\sup_{n} \frac{c_{n}}{b_{n}+c_{n}}A_{n}^{D}f(x)=D^{\ast}f(x)
\end{eqnarray*} 
Let \\
\begin{eqnarray*}
A=\{x:D^{\ast}f(x)\geq \lambda N_{\Phi}\}
\end{eqnarray*} \\
where $\ds N_{\Phi}$ will in this instance denote the Orlicz norm of $\ds f$.\\
\begin{eqnarray*}
\frac{c_{n}}{b_{n}+c_{n}}\left|A_{n}^{D}f(x)\right| & \leq & \frac{1}{b_{n}+c_{n}} \sum_{u\in \bigcup_{i=1}^{k-1}D_{i}\cap [0,n]}f(T^{u}x)\\
& \leq & \frac{1}{b_{n}+c_{n}} \sum_{u\in \bigcap_{i=1}^{k-1}D_{i}}f(T^{u}x)\\ 
& \leq & \frac{1}{l_{1}+\cdots+l_{k-1}}\sum_{u\in \bigcup_{i=1}^{k-1}}f(T^{u}x)\\
& =& R_{k-1}f(x)
\end{eqnarray*} 
Let $\ds T^{\ast}f(x)=\sup_{k}R_{k}f(x)$ and $\ds A_{k}=\{x:R_{k}f(x) \geq \lambda\}$, therefore\\
\begin{eqnarray*}
\mu(A^{\lambda}) \leq \sum_{k=1}^{\infty} \mu(A_{k}^{\lambda})
\end{eqnarray*}
Now if,\\
\begin{eqnarray*}
I & = & \int_{\{R_{k}f(x)\geq \lambda N_{\Phi}\}=A_{k}} \left(\frac{1}{\Phi(d_{1}+\cdots+d_{k})N_{\Phi}}\right)\sum_{u\in \bigcup_{i=1}^{k-1}D_{i}}f(T^{u}x)\mbox{d}x\\
& = & \int_{\{x:\sum_{u\in \bigcup_{i=1}^{k}D_{i}} T^{u}f(x) \geq \lambda N_{\Phi}(l_{1}+\cdots+l_{k})\}} \Phi\left(\frac{1}{N_{\Phi}(d_{1}+\cdots+d_{k})}\sum_{u \in \bigcup_{i=1}^{k}D_{i}} T^{u}f(x)\right)\mbox{d}x 
\end{eqnarray*}  
We have,\\
\begin{eqnarray*}
\mu(A_{k}^{\lambda})\Phi(\frac{\lambda (l_{1}+\cdots+l_{k})}{d_{1}+\cdots+d_{k}})\leq I \leq 1 \end{eqnarray*}
since \\
\begin{eqnarray*}
\left\|\sum_{u\in \bigcup_{i=1}^{k}D_{i}}T^{u}f(x)\right\|_{\Phi} \leq \frac{1}{N_{\Phi}(d_{1}+\cdots+d_{k})}
\end{eqnarray*}
Hence we have \\
\begin{eqnarray*}
\mu(A_{k}^{\lambda})\leq \frac{1}{\Phi(\frac{\lambda(l_{1}+\cdots+l_{k})}{d_{1}+\cdots+d_{k}})}
\end{eqnarray*}
Therefore if we have \\
\begin{eqnarray*}
\mu(\lambda) \leq \sum_{k=1}^{\infty} \mu(A_{k}^{\lambda}) = F(\lambda)
\end{eqnarray*} \\
For large enough $\ds \lambda$ we have, 
\begin{eqnarray*}
\frac{1}{\Phi(\frac{\lambda(l_{1}+\cdots + l_{k})}{d_{1}+\cdots+d_{k}})} \leq \frac{1}{ \Phi\left(\frac{l_{1}+\cdots+l_{k}}{d_{1}+\cdots+d_{k}}\right)} 
\end{eqnarray*}
Also,\\
\begin{eqnarray*}
\frac{1}{\Phi\left(\frac{\lambda(l_{1}+\cdots+l_{k})}{d_{1}+\cdots+d_{k}}\right)}\rightarrow 0
\end{eqnarray*}
monotonically as $\ds \lambda \rightarrow \infty$ for every $k$.\\
Therefore by the Lebesgue dominated convergence theorem $\ds F(\lambda)$ is an eventually monotone decreasing function that goes to $0$ as $\ds \lambda\rightarrow\infty$. \\
Since the maximal operator satisfies a weak-maximal inequality $\ds A_{n}^{D}f(x)$ \\
converges a.e.
\end{proof}
\newtheorem{second}[reinhold]{Proposition} 
\begin{second}
Let $\ds B_{k}$ and $\ds D_{k}$ as above.
Let $\ds l_{k}=\left|B_{k}\right|$ and $\ds d_{k}=|D_{k}|$.\\
Suppose that $\ds\forall k$ 
\begin{eqnarray*}
l_{1}+\cdots+l_{k} & \leq & Cl_{k+1}  \\
d_{k} & = & c_{k}l_{k}
\end{eqnarray*}
are such that $\ds \sum_{k=1}^{\infty} \frac{1}{\Phi\left(\frac{l_{k+1}}{l_{k}}\right)}\leq \infty$ and $\ds \sum_{k=1}^{\infty} \frac{1}{\Phi(\frac{1}{c_{k}})} \leq \infty$.\\
Then if $\ds\bigcup_{k=1}^{\infty} B_{k}$ is universally good  in $\ds\lp_{\Phi}$ then $\ds\bigcup_{k=1}^{\infty} (B_{k} \cup D_{k})$ is universally good in $\ds\lp_{\Phi}$.
\end{second}
\begin{proof}Choose $k_{0}$ so that $c_{k} \leq 1$ for all $k\geq k_{0}$.  Then
\begin{eqnarray*}
\frac{d_{1}+\cdots+d_{k}}{l_{1}+\cdots+l_{k}} & \leq & \frac{d_{1}+\cdots+d_{k_{0} -1}}{l_{1}+\cdots+l_{k}} + \frac{d_{k_{0}}+\cdots+d_{k-2}}{l_{1}+\cdots+l_{k}} + \frac{d_{k-1}+d_{k}}{l_{1}+\cdots+l_{k}} \\
& \leq & \frac{C_{0}}{l_{k}} + \frac{Cl_{k-1}}{l_{k}} +c_{k-1} + c_{k}
\end{eqnarray*} 
Therefore,
\begin{eqnarray*}
\frac{d_{1}+\cdots+d_{k}}{l_{1}+\cdots+l_{k}} & \leq & \frac{C_{0}}{l_{k}} + \frac{Cl_{k-1}}{l_{k}} + c_{k-1} +c_{k} ,\, \mbox{or}\\
\frac{l_{1}+\cdots+l_{k}}{d_{1}+\cdots+d_{k}} & \geq & \frac{1}{\frac{C_{0}}{l_{k}} + \frac{Cl_{k-1}}{l_{k}} + c_{k-1} + c_{k}} \\ 
& \geq & \frac{1}{4\max \left(\frac{C_{0}}{l_{k}},\frac{Cl_{k-1}}{l_{k}},c_{k-1},c_{k}\right)}= \min \left(\frac{1}{\frac{C_{0}}{l_{k}}},\frac{1}{\frac{Cl_{k-1}}{l_{k}}},\frac{1}{c_{k-1}},\frac{1}{c_{k}}\right)
\end{eqnarray*}  
Therefore,
\begin{eqnarray*}
\Phi\left(\frac{l_{1}+\cdots+l_{k}}{d_{1}+\cdots+d_{k}}\right) \geq \Phi\left(\frac{1}{4} \min(A_{k},B_{k},C{k},D_{k}\right)
\end{eqnarray*}  or,
\begin{eqnarray*}
\frac{1}{\Phi\left(\frac{l_{1}+\cdots+l_{k}}{d_{1}+\cdots+d_{k}}\right)} & \leq &  \frac{1}{\Phi\left(\frac{1}{4} \min(A_{k},B_{k},C_{k},D_{k})\right)} \\
& =&  \max\left(\frac{1}{\Phi(\frac{1}{4}A_{k})},\frac{1}{\Phi(\frac{1}{4}B_{k})},\frac{1}{\Phi(\frac{1}{4}B_{k})},\frac{1}{\Phi(\frac{1}{4}D_{k})}\right)
\end{eqnarray*}
and 
\begin{eqnarray*}
\sum_{k=1}^{\infty} \frac{1}{\Phi(\frac{l_{1}+\cdots+l_{k}}{d_{1}+\cdots+d_{k}})}< \sum_{k=1}^{\infty} \max (P_{k},Q_{k},R_{k},S_{k})<\infty
\end{eqnarray*}
\end{proof}
\section{}
Suppose $f$ is a monotone decreasing function on $(0,1)$. \\ 
Let 
\begin{eqnarray*}
B_\lambda=\{x:A_N=\frac{1}{N}\sum_{k=1}^Nf(T^{n_k}x)\geq\lambda\}.
\end{eqnarray*}
For each $k$ where $\ds 0\leq k \leq N$, let 
\begin{eqnarray*}
a_k = \sup\{x:y=T^{-n_k}(x)\in B_\lambda\: \mbox{and}\: T^{n_k}(x)\leq T^{n_p}(x)\quad \forall\quad 1\leq p \leq N \}.
\end{eqnarray*}
Intuitively this is the supremum of the $x$ values such that there is a $y$ with $\ds T^{n_k}(y)=x$ and $x$ is the smallest distance of the partial orbit $\ds\{T^{n_k}(y)\}_{k=1}^{N}$ to the origin.
Let $\ds S_k= T^{-n_k}([0,a_k]).$
\newtheorem{union}{Theorem}[section]
\begin{union}
\label{firstfourthsection}
$\ds \bigcup_{k=0}^{N-1} S_k=B_\lambda$ up to a set of measure zero, the union being disjoint.
\end{union}
\begin{proof}
Suppose $\ds x\in A_k$ for some $\ds 1\leq k\leq N$. 
Then $\ds T^{n_k}(x)\geq  T^{n_p}(x)$ $\forall 1\leq p\leq N$. 
Therefore the $\ds \supseteq$ inclusion has been proved.
Now suppose $\ds x\in B$. 
There is a point of the set $\ds\{T^{n_k}(x)\}_{k=1}^{N}$ that is closest to the origin, say $\ds T^{n_k}(x)$. If $\ds T^{n_p}(x)>a_p$, we contradict the definition
 of $a_p$.
Thus $x\in A_p$.  
It remains to prove the disjointness assertion. To this end suppose that 
$\ds A_p\bigcap A_q \neq \emptyset$ and that $\ds a_p>a_q$.  
At this point it may be convient to view modulo $1$ arithmetic on $[0,1)$ as a rotation of the circle.  See diagrams below. 
\newpage
\begin{figure}
\begin{center}
\scalebox{1} 
{
\begin{pspicture}(0,-2.328125)(10.622812,2.328125)
\pscircle[linewidth=0.04,dimen=outer](1.7609375,0.1496875){0.72}
\pscircle[linewidth=0.04,dimen=outer](4.1409373,1.1696875){0.74}
\pscircle[linewidth=0.04,dimen=outer](3.8609376,-1.1303124){0.72}
\psline[linewidth=0.04cm,arrowsize=0.05291667cm 2.0,arrowlength=1.4,arrowinset=0.4]{->}(2.6609375,0.4496875)(3.2009375,0.7896875)
\psline[linewidth=0.04cm,arrowsize=0.05291667cm 2.0,arrowlength=1.4,arrowinset=0.4]{->}(2.6009376,-0.3303125)(3.0609374,-0.6703125)
\psdots[dotsize=0.12](2.4809375,0.1496875)
\psdots[dotsize=0.12](2.1209376,0.7696875)
\psdots[dotsize=0.12](1.3609375,0.7296875)
\psdots[dotsize=0.12](1.0609375,0.1296875)
\psdots[dotsize=0.12](1.3209375,-0.3903125)
\psdots[dotsize=0.12](4.8809376,1.2496876)
\psdots[dotsize=0.12](4.5809374,1.7496876)
\psdots[dotsize=0.12](4.0609374,1.9096875)
\psdots[dotsize=0.12](3.5809374,1.6096874)
\psdots[dotsize=0.12](4.5609374,-1.0903125)
\psdots[dotsize=0.12](4.4409375,-1.4303125)
\psdots[dotsize=0.12](4.2609377,-0.5103125)
\psdots[dotsize=0.12](3.8209374,-0.4103125)
\usefont{T1}{ptm}{m}{n}
\rput(5.6423435,1.2596875){$\ds e_p$ at $0$}
\usefont{T1}{ptm}{m}{n}
\rput(4.702344,1.7996875){$\ds e_q$}
\usefont{T1}{ptm}{m}{n}
\rput(4.092344,2.1396875){$\ds y_q$}
\usefont{T1}{ptm}{m}{it}
\rput(3.3314064,1.8196875){$\ds y_p$}
\usefont{T1}{ptm}{m}{n}
\rput(5.282344,-1.1003125){$\ds e_q$ at $0$}
\usefont{T1}{ptm}{m}{n}
\rput(4.7123437,-0.5603125){$\ds y_q$}
\usefont{T1}{ptm}{m}{n}
\rput(4.5823436,-1.5603125){$\ds e_p$}
\usefont{T1}{ptm}{m}{n}
\rput(3.6123438,-0.2803125){$\ds y_p$}
\usefont{T1}{ptm}{m}{n}
\rput(2.4723437,-0.8203125){$\ds T^q$}
\usefont{T1}{ptm}{m}{n}
\rput(2.6723437,0.9796875){$\ds T^p$}
\usefont{T1}{ptm}{m}{n}
\rput(2.6523438,0.1196875){$0$}
\usefont{T1}{ptm}{m}{n}
\rput(2.0223436,0.9196875){$\ds e_p$}
\usefont{T1}{ptm}{m}{n}
\rput(1.1223438,0.8596875){$\ds e_q$}
\usefont{T1}{ptm}{m}{n}
\rput(0.67234373,0.2796875){$\ds y_q$}
\usefont{T1}{ptm}{m}{n}
\rput(0.9323437,-0.5803125){$\ds y_p$}
\usefont{T1}{ptm}{m}{n}
\rput(5.7523437,-2.1003125){$\ds y_p$ is brought closer to $0$ by $\ds T^q$ than by $\ds T^p$}
\end{pspicture} 
}
\end{center}
\caption{Theorem ~\ref{firstfourthsection}}
\end{figure}

\begin{figure}
\begin{center}
\scalebox{1} 
{
\begin{pspicture}(0,-2.328125)(10.622812,2.328125)
\pscircle[linewidth=0.04,dimen=outer](1.7609375,0.1496875){0.72}
\pscircle[linewidth=0.04,dimen=outer](4.1409373,1.1696875){0.74}
\pscircle[linewidth=0.04,dimen=outer](3.8609376,-1.1303124){0.72}
\psline[linewidth=0.04cm,arrowsize=0.05291667cm 2.0,arrowlength=1.4,arrowinset=0.4]{->}(2.6609375,0.4496875)(3.2009375,0.7896875)
\psline[linewidth=0.04cm,arrowsize=0.05291667cm 2.0,arrowlength=1.4,arrowinset=0.4]{->}(2.6009376,-0.3303125)(3.0609374,-0.6703125)
\psdots[dotsize=0.12](2.4809375,0.1496875)
\psdots[dotsize=0.12](2.1209376,0.7696875)
\psdots[dotsize=0.12](1.3609375,0.7296875)
\psdots[dotsize=0.12](1.0609375,0.1296875)
\psdots[dotsize=0.12](1.3209375,-0.3903125)
\psdots[dotsize=0.12](4.8809376,1.2496876)
\psdots[dotsize=0.12](4.5809374,1.7496876)
\psdots[dotsize=0.12](4.0609374,1.9096875)
\psdots[dotsize=0.12](3.5809374,1.6096874)
\psdots[dotsize=0.12](4.5609374,-1.0903125)
\psdots[dotsize=0.12](4.4409375,-1.4303125)
\psdots[dotsize=0.12](4.2609377,-0.5103125)
\psdots[dotsize=0.12](3.8209374,-0.4103125)
\usefont{T1}{ptm}{m}{n}
\rput(5.6423435,1.2596875){$\ds e_p$ at $0$}
\usefont{T1}{ptm}{m}{n}
\rput(4.702344,1.7996875){$\ds e_q$}
\usefont{T1}{ptm}{m}{n}
\rput(4.092344,2.1396875){$\ds y_p$}
\usefont{T1}{ptm}{m}{it}
\rput(3.3314064,1.8196875){$\ds y_q$}
\usefont{T1}{ptm}{m}{n}
\rput(5.282344,-1.1003125){$\ds e_q$ at $0$}
\usefont{T1}{ptm}{m}{n}
\rput(4.7123437,-0.5603125){$\ds y_p$}
\usefont{T1}{ptm}{m}{n}
\rput(4.5823436,-1.5603125){$\ds e_p$}
\usefont{T1}{ptm}{m}{n}
\rput(3.6123438,-0.2803125){$\ds y_q$}
\usefont{T1}{ptm}{m}{n}
\rput(2.4723437,-0.8203125){$\ds T^q$}
\usefont{T1}{ptm}{m}{n}
\rput(2.6723437,0.9796875){$\ds T^p$}
\usefont{T1}{ptm}{m}{n}
\rput(2.6523438,0.1196875){$0$}
\usefont{T1}{ptm}{m}{n}
\rput(2.0223436,0.9196875){$\ds e_p$}
\usefont{T1}{ptm}{m}{n}
\rput(1.1223438,0.8596875){$\ds e_q$}
\usefont{T1}{ptm}{m}{n}
\rput(0.67234373,0.2796875){$\ds y_p$}
\usefont{T1}{ptm}{m}{n}
\rput(0.9323437,-0.5803125){$\ds y_q$}
\usefont{T1}{ptm}{m}{n}
\rput(5.7523437,-2.1003125){$\ds y_p$ is brought closer to $0$ by $\ds T^q$ than by $\ds T^p$}
\end{pspicture} 
}
\end{center}
\caption{Theorem ~\ref{firstfourthsection}}
\end{figure}

Note that the rotation of the circle is orientation preserving.  From these diagrams it is clear that an intersection of these sets must result in a contradiction of the definition of either $a_p$ or $a_q$.\\
Let $L$ denote the measure of $\ds B_\lambda$.  We have that $\ds B_\lambda=\bigcup_{k=1}^{N} I_k$, where $\ds I_k$ is an \\
interval, possibly empty, and if $\ds [c_k,d_k]$ denotes such an interval then $\ds T^{n_k}(c_k)=0$ and $\ds T^{n_k}(d_k)=a_k$ where $a_k$ is as above.  \\
We now create an interval of length $L$ which consists of intervals $\{J_k\}_{k=1}^{N}$ with $\left| {J_k}\right|=\left|{I_k}\right|$, and such that the orientation of the $\{J_k\}$ is the same as that of the ${I_k}$.  See diagram below.  Let us call this new space $X$.
			Map $B_\lambda$ to $X$ as follows.
			Let $\ds\Phi :B_\lambda  \rightarrow X$ where $\ds\Phi (I_k) = J_k $, where $\ds\Phi$ is defined in the obvious way as an orientation preserving isometry when so restricted.\\
	We now define a sequence of measure preserving transformations $\ds\{ \Psi _k \} _{k=1}^{N}$  on the probability space $\ds(X,\mathcal{B},\frac{\mu}{L})$, where $\mu$ is the lebesgue measure of the unit interval.\\
If $\ds J_k=[r_k,s_k]$ we let $\ds\Psi_k(x) = x+(L-r_k)\mbox{mod}L$, so that $\ds\Psi ( r_k) = 0$ and $\ds\Psi (s_k) =a_k$.   Now let $\ds F_N(x) =\frac{1}{N} \sum _{k=1}^{N} f(\Psi _k (x))$ for x in $\ds[0,L]$.
\end{proof} 
\newtheorem{prop}[union]{Porposition}
 \begin{prop}
 \label{secondfourthsection}
 Let $C=\{x\in [0,L] : F_N(x)\geq  \lambda \}$.  Then $\left| C \right| = L$.\\     
\end{prop}
\begin{proof}
Let $x\in J_k$.  Since $\ds x\in J_k $ we have that  $\ds y=\Phi^{-1}(x)\in I_k$. \\
It is true that $\ds\forall 1\:\leq p\leq N$, we have $\ds T^{n_p}(y)\geq \Psi ^{p}(x),$ and \\
therefore by the monotonicity of $f$ $\ds f(T^{n_k})(y)\leq f(\Psi ^{k}(x))$ and hence also $F_N(x)\geq A_N(y)\geq\lambda$.  \\
The first assertion of the last line follows from the fact that the transformations $T^{n_p}$ and $\Psi_k$ map $I_p$ and $J_k$ to the interval $[0,a_p]$ respectively and therefore the will be the same number of the intervals from the collections, however in $X$ we have eliminated the space between the intervals and thus the distance from each point to the origin has been decreased.\\
\begin{figure}[h]
\begin{center}
\scalebox{1} 
{
\begin{pspicture}(0,-1.3892188)(8.122812,1.4292188)
\pscircle[linewidth=0.04,dimen=outer](1.4909375,-0.09921875){1.29}
\pscircle[linewidth=0.04,dimen=outer](6.3009377,-0.18921874){0.88}
\psline[linewidth=0.04cm,arrowsize=0.05291667cm 2.0,arrowlength=1.4,arrowinset=0.4]{->}(3.2409375,-0.12921876)(5.0209374,-0.10921875)
\psline[linewidth=0.04cm](0.5209375,0.97078127)(0.6809375,0.79078126)
\psline[linewidth=0.04cm](0.9209375,1.1307813)(1.0209374,0.97078127)
\psline[linewidth=0.04cm](2.1209376,1.1507813)(2.0009375,0.97078127)
\psline[linewidth=0.04cm](2.2809374,0.81078124)(2.4009376,0.97078127)
\psline[linewidth=0.04cm](2.5209374,0.43078125)(2.7609375,0.55078125)
\psline[linewidth=0.04cm](2.6609375,0.19078125)(2.8409376,0.27078125)
\psline[linewidth=0.04cm](2.6409376,-0.00921875)(2.9009376,0.07078125)
\psline[linewidth=0.04cm](2.6609375,-0.38921875)(2.8809376,-0.46921876)
\psline[linewidth=0.04cm](7.0609374,-0.18921874)(7.3009377,-0.18921874)
\psline[linewidth=0.04cm](7.0209374,0.09078125)(7.1809373,0.25078124)
\psline[linewidth=0.04cm](6.8209376,0.35078126)(6.9609375,0.57078123)
\usefont{T1}{ptm}{m}{n}
\rput(2.8923438,-0.17921875){$0$}
\usefont{T1}{ptm}{m}{n}
\rput(3.0523438,-0.35921875){$\ds I_r$}
\usefont{T1}{ptm}{m}{n}
\rput(2.8823438,0.46078125){$\ds I_s$}
\usefont{T1}{ptm}{m}{n}
\rput(2.1823437,1.1807812){$\ds I_1$}
\usefont{T1}{ptm}{m}{n}
\rput(0.60234374,1.2407813){$\ds I_t$}
\usefont{T1}{ptm}{m}{n}
\rput(4.2323437,-0.51921874){$\ds \Phi$}
\usefont{T1}{ptm}{m}{n}
\rput(7.3923435,-0.01921875){$\ds J_1$}
\usefont{T1}{ptm}{m}{n}
\rput(7.072344,0.44078124){$\ds J_s$}
\psdots[dotsize=0.12](2.7809374,-0.20921876)
\end{pspicture} 
}
\end{center}
\caption{Proposition ~\ref{secondfourthsection}}
\end{figure}

\newpage

\begin{figure}[h]
\begin{center}
\scalebox{1} 
{
\begin{pspicture}(0,-1.2289063)(6.2228127,1.2089063)
\pscircle[linewidth=0.04,dimen=outer](1.8409375,0.00890625){0.96}
\pscircle[linewidth=0.04,dimen=outer](4.4009376,0.00890625){0.96}
\psline[linewidth=0.04cm](2.4609375,0.80890626)(2.3809376,0.74890625)
\psline[linewidth=0.04cm](2.6409376,0.5289062)(2.7009375,0.5889062)
\psline[linewidth=0.04cm](0.9809375,-0.49109375)(1.1009375,-0.43109375)
\psline[linewidth=0.04cm](1.2609375,-0.79109377)(1.3609375,-0.69109374)
\psline[linewidth=0.04cm](1.6609375,-0.79109377)(1.6609375,-0.93109375)
\psline[linewidth=0.04cm](5.2609377,0.00890625)(5.4809375,0.00890625)
\psline[linewidth=0.04cm](5.2209377,0.32890624)(5.3609376,0.38890624)
\psline[linewidth=0.04cm](4.9809375,0.64890623)(5.1009374,0.76890624)
\psdots[dotsize=0.12](1.8609375,1.1289062)
\psdots[dotsize=0.12](1.2009375,0.9089062)
\psdots[dotsize=0.12](0.7809375,0.28890625)
\psline[linewidth=0.04cm](2.6409376,0.04890625)(2.8809376,0.04890625)
\usefont{T1}{ptm}{m}{n}
\rput(2.6123438,0.8389062){$\ds J_s$}
\usefont{T1}{ptm}{m}{n}
\rput(2.8923438,0.39890626){$\ds J_1$}
\usefont{T1}{ptm}{m}{n}
\rput(0.65234375,-0.60109377){$\ds J_p$}
\usefont{T1}{ptm}{m}{n}
\rput(1.0623437,-1.0010937){$\ds J_r$}
\usefont{T1}{ptm}{m}{n}
\rput(3.1323438,-0.36109376){$\ds \Psi_p$}
\usefont{T1}{ptm}{m}{n}
\rput(5.492344,0.17890625){$\ds J_p$}
\usefont{T1}{ptm}{m}{n}
\rput(5.382344,0.67890626){$\ds J_r$}
\psline[linewidth=0.04cm,arrowsize=0.05291667cm 2.0,arrowlength=1.4,arrowinset=0.4]{->}(2.9209375,-0.09109375)(3.3609376,-0.09109375)
\end{pspicture} 
}
\end{center}
\caption{Proposition ~\ref{secondfourthsection}}
\end{figure}

\end{proof}

\newtheorem{lam}[union]{Theorem}
\begin{lam}
Let $\ds M_\lambda = \sup \{L : \frac{1}{L} \int_{0}^L f(x)\mbox{d}x > \lambda \}.$ Then $\ds\left|B_\lambda \right|\leq M_\lambda.$
\end{lam}
\begin{proof}
Since $\ds\left|B_\lambda \right|= \left| C\right| = L$ , as above and $\ds\Psi _{k}$ is an m.p.t. of the space $X$, we have by the change of variable formula\\
\begin{eqnarray*}
\lambda \leq  \frac{1}{L} \int_{0}^{L} F(x)\mbox{d}x & = & \frac{1}{L} \int_{0}^{L} \sum_{k=0}^{N-1}f(\Psi_{k}(x))\mbox{d}x \\
& = & \frac{1}{L} \int_{0}^{L} f(x) \mbox{d}x.
\end{eqnarray*}
Therefore $\ds \lambda\leq M_\lambda$.
\end{proof}
\newtheorem{lam2}[union]{Theorem}
\begin{lam2}
\label{fourthfourthsection}
Let $\ds P= \sup_ {n_{0}<n_{1}<\cdots<n_{k-1}} \left| B_\frac{\lambda}{2}\right|$. Then $\ds P= M_{\lambda}.$  
\end{lam2}
\begin{proof}
Let $\eta,\delta>0$. Since $f$ is a monotone decreasing function $\ds f\chi _{[\epsilon, M_{\lambda}]}$  is a bounded funcction and therefore Riemann Integrable.
Therefore there exists a number $\ds r_{k}$ such that if $\ds[\epsilon,M_{\lambda}]$  is partitioned into $\ds r_{k}$ intervals of equal lenghth $\ds\{ I_{j}\}_{j=1}^{r_{k}}$
we have $$\ds\left| \sum_{i=1}^{r_{k}} f(x_{i}) \left| I_{i}\right|- \int_{\epsilon}^ {M_{\lambda}} f(x) \mbox{d}x \right|<\eta.$$
Now let $\ds I_{j}=[a_{j},b_{j}]$, $\ds I_{1}$ has right endpoint $\ds M_{\lambda}$ and $\ds I_{r_{k}}$ has left endpont $\epsilon$.\\
Choose $\ds nf(x)dx_{j}$  so that $\ds T^{n_{j}}(M_{\lambda})\in I_{j}$ and it's distance from $b_{j}$ is less than some number $\beta >0$ for all $\ds 1\leq j \leq r_{k}$.  Now choose $n_{j}$,  $\ds r_{j}\leq j \leq 2r_{j}$ so that $\ds T^{n_{j}}(\epsilon)\in I_{j}$ and is within some $\beta$ of $a_{j}.$  For all $x$ except for those which are contained in a set whose measure is determined by $\beta$ we have for $\ds x\in I_{j}$,  $\ds T^{n_{i}}(x)\in I_{j+i}$ for $\ds 1\leq i\leq r_{k}-j$. Also $\ds T^{n_{i}}(x)\in I_{j-i}$ for $r_{k}\leq i \leq r_{k} +j$.\\
See diagram.  
\begin{figure}[h]
\begin{center}
\scalebox{1} 
{
\begin{pspicture}(0,-2.438125)(8.13,2.438125)
\psline[linewidth=0.04cm](0.6809375,1.8996875)(6.6809373,1.8996875)
\psline[linewidth=0.04cm](0.7409375,-0.1003125)(6.7609377,-0.1203125)
\psline[linewidth=0.04cm](0.8409375,-1.9203125)(6.8809376,-1.9003125)
\psline[linewidth=0.04cm](0.7009375,2.0796876)(0.7009375,1.8196875)
\psline[linewidth=0.04cm](0.7409375,0.0)(0.7409375,-0.1803125)
\psline[linewidth=0.04cm](0.8409375,-1.7603126)(0.8409375,-2.0203125)
\psline[linewidth=0.04cm](6.7009373,2.0996876)(6.7009373,1.7396874)
\psline[linewidth=0.04cm](6.7409377,0.0196875)(6.7409377,-0.2603125)
\psline[linewidth=0.04cm](6.8609376,-1.7603126)(6.8609376,-2.0803125)
\psline[linewidth=0.04cm](6.2809377,1.9996876)(6.2809377,1.8396875)
\psline[linewidth=0.04cm](6.2609377,-0.0203125)(6.2609377,-0.2003125)
\psline[linewidth=0.04cm](5.8209376,-0.0403125)(5.8409376,-0.1603125)
\psline[linewidth=0.04cm](5.4009376,-0.0403125)(5.4209375,-0.2203125)
\psline[linewidth=0.04cm](4.2609377,0.0796875)(4.2609377,-0.2203125)
\psline[linewidth=0.04cm](3.8409376,0.0596875)(3.8409376,-0.1803125)
\psline[linewidth=0.04cm](3.4809375,0.0396875)(3.4809375,-0.2003125)
\psline[linewidth=0.04cm](3.0409374,0.0396875)(3.0409374,-0.2203125)
\psline[linewidth=0.04cm](3.9609375,2.0196874)(3.9409375,1.7996875)
\psline[linewidth=0.04cm](3.7009375,2.0596876)(3.7009375,1.8196875)
\psline[linewidth=0.04cm](1.1009375,2.0596876)(1.1009375,1.8196875)
\psline[linewidth=0.04cm](1.3809375,2.0596876)(1.3609375,1.8596874)
\psline[linewidth=0.04cm](1.6809375,2.0396874)(1.6809375,1.7996875)
\psline[linewidth=0.04cm](1.9809375,2.0196874)(1.9809375,1.7996875)
\psline[linewidth=0.04cm](0.9609375,0.0)(0.9609375,-0.1803125)
\psline[linewidth=0.04cm](1.2009375,0.0196875)(1.2009375,-0.2003125)
\psline[linewidth=0.04cm](1.4009376,0.0796875)(1.4009376,-0.2003125)
\psline[linewidth=0.04cm](1.6009375,0.0596875)(1.6009375,-0.1803125)
\psline[linewidth=0.04cm](1.1809375,-1.7803125)(1.1809375,-2.0203125)
\psline[linewidth=0.04cm](1.3809375,-1.8003125)(1.3809375,-2.0203125)
\psline[linewidth=0.04cm](1.7009375,-1.8203125)(1.7009375,-2.0203125)
\psline[linewidth=0.04cm](2.0009375,-1.7603126)(2.0009375,-2.0603125)
\psline[linewidth=0.04cm](3.6209376,-1.7603126)(3.6209376,-2.0403125)
\psline[linewidth=0.04cm](3.9009376,-1.8003125)(3.9009376,-2.0403125)
\psline[linewidth=0.04cm](4.2409377,-1.8003125)(4.2409377,-2.0803125)
\psline[linewidth=0.04cm](4.5209374,-1.8003125)(4.5209374,-2.0203125)
\psarc[linewidth=0.04](4.2709374,-1.9503125){0.17}{0.0}{180.0}
\psarc[linewidth=0.04](3.9209375,-1.9003125){0.18}{0.0}{180.0}
\psarc[linewidth=0.04](1.6909375,-1.8903126){0.17}{0.0}{180.0}
\psarc[linewidth=0.04](1.3509375,-1.8303125){0.15}{0.0}{180.0}
\psarc[linewidth=0.04](3.8909376,-0.1303125){0.31}{0.0}{180.0}
\psarc[linewidth=0.04](3.3809376,-0.0803125){0.18}{0.0}{180.0}
\psarc[linewidth=0.04](6.3409376,-0.1603125){0.36}{0.0}{180.0}
\psarc[linewidth=0.04](5.7409377,-0.1803125){0.24}{0.0}{180.0}
\psarc[linewidth=0.04,linestyle=dotted,dotsep=0.16cm](4.7509375,0.1896875){0.47}{0.0}{180.0}
\psline[linewidth=0.04cm,arrowsize=0.05291667cm 2.0,arrowlength=1.4,arrowinset=0.4]{->}(5.7209377,1.0996875)(4.9409375,1.0996875)
\usefont{T1}{ptm}{m}{n}
\rput(6.7559376,1.5696875){$\ds M_{\lambda}$}
\usefont{T1}{ptm}{m}{n}
\rput(6.4223437,2.2296875){$\ds I_1$}
\usefont{T1}{ptm}{m}{n}
\rput(3.7423437,2.2296875){$\ds I_j$}
\usefont{T1}{ptm}{m}{n}
\rput(1.3659375,2.2496874){$\ds I_{r_k}$}
\usefont{T1}{ptm}{m}{n}
\rput(0.95234376,1.6696875){$\ds \epsilon$}
\usefont{T1}{ptm}{m}{n}
\rput(1.0123438,-0.3703125){$\ds \epsilon$}
\usefont{T1}{ptm}{m}{n}
\rput(1.2123437,-2.1703124){$\ds \epsilon$}
\usefont{T1}{ptm}{m}{n}
\rput(1.1323438,-1.5303125){$\ds I_{r_k}$}
\usefont{T1}{ptm}{m}{n}
\rput(1.7359375,-1.5103126){$\ds I_{r_{k-1}}$}
\usefont{T1}{ptm}{m}{n}
\rput(2.0159376,-2.2103126){$\ds I_{r_{k-2}}$}
\usefont{T1}{ptm}{m}{n}
\rput(6.4023438,-0.3303125){$\ds I_1$}
\usefont{T1}{ptm}{m}{n}
\rput(5.862344,-0.4503125){$\ds I_2$}
\usefont{T1}{ptm}{m}{n}
\rput(5.4223437,-0.2503125){$\ds I_3$}
\psline[linewidth=0.04cm,arrowsize=0.05291667cm 2.0,arrowlength=1.4,arrowinset=0.4]{->}(2.2809374,-1.1603125)(3.3609376,-1.1803125)
\end{pspicture} 
}
\end{center}
\caption{Theorem ~\ref{fourthfourthsection}}
\end{figure}

Therefore with $\ds x_{i}\in I_{i}$,
\begin{eqnarray*}
\frac{1}{2r_{k}} \sum_{i=1}^{2r_{k}} f(T^{n_{i}}(x)) & \geq &  \frac{1}{2r_{k}} \sum_{i=r_{k}}^{r_{k}+j} f(T^{n_{i}}(x)) + \frac{1}{2r_{k}} \sum_{i=1}^{r_{k}-j} f(T^{n_{i}}(x)) \\
& = & \frac{1}{2r_{k}} \sum_{i=1}^{r_{k}} f(x_{i}) \\
& = & \frac{1}{2} \frac{1}{M_{\lambda} - \epsilon}\frac{M_{\lambda} - \epsilon}{r_{k}} \sum_{i=1}^{r_{k}} f(x_{i}) \\
& = & \frac{1}{2} \frac{1}{M_\lambda - \epsilon} \sum_{i=1}^{r_{k}} f(x_{i})\left|I_{i}\right| = F(x,\epsilon).
\end{eqnarray*}
\begin{eqnarray*}
\left|2F(x,\epsilon) - \frac{1}{M_{\lambda}}\int_{0}^{M_{\lambda}} f(x)\mbox{d}x\right| & \leq & \left|2F(x,\epsilon)-\frac{1}{M_{\lambda}}\int_{\epsilon}^{M_{\lambda}}f(x)\mbox{d}x\right|+\frac{1}{M_{\lambda}}\int_{0}^{\epsilon}f(x)\mbox{d}x\\
& < & \left|2F(x,\epsilon)-\frac{1}{M_{\lambda}}\int_{\epsilon}^{M_{\lambda}}f(x)dx\right|+\frac{\eta}{M_{\lambda}}.
\end{eqnarray*}
\\Now,	
\begin{eqnarray*}
\left|2F(x,\epsilon)-\frac{1}{M_{\lambda}}\int_{\epsilon}^{M_{\lambda}}f(x)\mbox{d}x\right|& \leq & \left|2F(x,\epsilon)-\frac{1}{M_{\lambda}- \epsilon}\int_{\epsilon}^{M_{\lambda}}f(x)\mbox{d}x\right|\\
& + & \left|\frac{1}{M_{\lambda}-\epsilon}\int_{\epsilon}^{M_{\lambda}}f(x)\mbox{d}x  \frac{1}{M_{\lambda}}\int_{\epsilon}^{M_{\lambda}}f(x)\mbox{d}x\right| \\
& \leq & \eta +\left|\frac{M_{\lambda}-(M_{\lambda}-\epsilon)}{(M_{\lambda}-\epsilon)M_{\lambda}}\right|\|f\|_1\\
& < & \eta + \frac{\epsilon}{M_{\lambda}}\|f\|_1.
\end{eqnarray*}
Therefore,\\
\begin{eqnarray*}
2F(x,\epsilon)> \frac{1}{M_{\lambda}}\int_{0}^{M_{\lambda}}f(x)dx - \frac{\eta}{M_{\lambda}}-\eta - \frac{\epsilon}{M_{\lambda}}\|f\|_1. 
\end{eqnarray*}
Choosing $\epsilon$ and $\eta$ small enough gives:\\
\begin{eqnarray*}
2\frac{1}{2r_{k}}\sum_{i=1}^{2r_{k}}f(T^{n_{i}}x)\geq 2F(x,\epsilon)>\frac{1}{M_\lambda}\int_{0}^{M_{\lambda}}f(x)\mbox{d}x>\lambda.
\end{eqnarray*}
This implies that for the finite subsequence $\ds n_{1}<\cdots<n_{2r_{k}}$ we have that 
\begin{eqnarray*}
 \left|\left\{x:\frac{1}{2r_{k}}\sum_{i=1}^{2r_{k}}f(T^{n_{i}}x)>\frac{\lambda}{2}\right\}\right|\geq M_{\lambda}-\delta,
 \end{eqnarray*}
where $\delta$ is arbitrarily small, assuming that $\epsilon$ and $\beta$ have been made sufficiently small.   
\end{proof}

\newtheorem{finite}[union]{Theorem}	  
\begin{finite}
Given any interval $I$ of length $\ds M_\lambda$ and any $\delta>0$  there exists a finite subsequence of integers $\ds n_{0}<n_{1}<\cdots<n_{k-1}$ and a subinterval $\ds I_{\delta} \subseteq I$, $\ds\left|I_{\delta} \right|>M_{\lambda} - \delta $ such that\\ $\forall x\in I_\delta$  
\begin{eqnarray*}
\frac{1}{k} \sum_{j=0}^{k-1} f(T^{n_{J}(x)}) \geq \frac{\lambda}{2}.
\end{eqnarray*}  
Furthermore the sequence can be made arbitrarily long. Also the choice of $n_{0}$ can be taken arbitrarily large.
\end{finite}
\begin{proof}
Let $\ds I=[a,a+M_{\lambda}]$.  In the previous proof replace $\epsilon$ by $ a+\epsilon$ and $M_{\lambda}$ by $a + M_{\lambda}$ and map these into the partitions of $[\epsilon,M_{\lambda}]$.  The ergodicity of the transformation ensures the claim regarding $n_{0}$.  By refining the partition one creates more intervals and the sequence can be made longer.
\end{proof}
\newtheorem{blocks}[union]{Theorem}
\begin{blocks}
Suppose that $f$ is a monotone decreasing function on $(0,1)$ and there exists a sequence $\ds{s_{k}}(s_{k}\rightarrow \infty)$ such that if 
\begin{eqnarray*}
a_{k}=\sup \{\lambda : \frac{1}{\lambda}\int_{0}^{\lambda}f(x)\mbox{d}x>\frac{s_{k}}{2}\}
\end{eqnarray*} 
we have $\ds \sum_{i=1}^{\infty}a_{k} = \infty$.\\
Then if ${c_{k}}$ is a sequence such that $\ds\frac{s_{k}}{c_{k}}\rightarrow \infty$ then there exists a block sequence $\ds bigcup B_{k}$ and a perturbation of this sequence $\ds\bigcup (B_{k}\cup D_{k})$ where $\ds\left|D_{k}\right|=\frac{1}{c_{k}}\left|B_{k}\right|$ such that the ergodic averages of $f$ along this subsequence fail to converge a.e.
\end{blocks}
\begin{proof}
	Let $\ds p_{k}=\sum_{j=1}^{k}a_{j}\mbox{mod}(1)$, and $\ds J_{K}=[p_{k},p_{k+1}]$. Since $\ds \sum{a_j}$ diverges each point of $[0,1)$ is in infinitely many of the $J_{k}$. We construct the sequence inductively as follows:\\
	let $\ds n_{1}=1$, $\ds l_{1}=1$, $\ds d_{1}=1$ and suppose that $\ds n_{1},\cdots,n_{k-1}$,$\ds l_{1},\cdots,l{k-1}$,
	$\ds D_{1},\cdots,D_{k-1}$,$\ds d_{1},\cdots,d_{k-1}$ have already been chosen. The block $B_{k}$ will satisfy the following
	\begin{enumerate}
	\item $\ds l_{k}>n_{k-1}$
	\item $\ds l_{k}\geq kl_{k-1}\geq l_{1}+\cdots+l_{k-1}.$
	\end{enumerate}
	Given $\ds\delta _{k}$ choose an integer $\ds d_{k}$ large enough so that there exists a subsequence of length $\ds d_{k}$ where
	\begin{eqnarray*}
	 \frac{1}{d_{k}} \sum_{0}^{d_{k}-1}f(T^{n_{j}}x)>\frac{s_{k}}{2}\quad \forall \: x\in (J_{k-1})_{\delta_{k}} . 
	 \end{eqnarray*}
	Now $\ds d_{k}$ and $\ds l_{k}$ may be chosen so that $\ds d_{k}=c_{k}l_{k}$ and the above conditions are satisfied.\\
	Note the fact that we may arbitrarily lenghthen a subsequence is key to finding the integer $d_{k}$.  Let $B_{k}$ consist of a block of integers starting to the right of $\ds D_{k-1}$ and $\ds D_{k}$ be $\ds d_{k}$ integers to the right of $B_{k}$ that yield the above inequality for the points in $\ds (J_{k-1})_{\delta_{k}}$.\\
	Therefore, $\ds\:\forall\: x\in (J_{k-1})_{\delta_{k}}$\\
	\begin{eqnarray*}
	\frac{1}{l_{1}+\cdots+l_{k}+d_{1}+\cdots+d_{k}}\sum_{u\in \bigcup (B_{j}\cup D_{j})}f(T^{u}x)
	& \geq & C \frac{1}{l_{k}}\sum_{u\in D_{k}}f(T^{u}x)\\
	& \geq & C \frac{d_{k}}{l_{k}}s_{k} \\
	& = &  Cs_{k}c_{k}\rightarrow \infty .
	\end{eqnarray*}
	 \\
  If the $\ds\delta _{k}$'s are chosen small enough, there will exist a set of positive measure $J$ so that each $x\in J$  is in infinitely many of the $\ds (J_{k-1})_{\delta _{k}}$.\\
  Clearly such a point will have a subsequence of averages which diverge to infinity.
  \end{proof}
  \theoremstyle{remark}
  \newtheorem{ex}{Example}[section]
  \begin{ex}
  Let $s>0$, define 
  \begin{eqnarray*}
  s_{k} & = & \frac{k}{(\log{k}) ^{s-1} (\log{\log{k}})^{s}}\quad  \mbox{and}\\
  c_{k} & = & \frac{\log\log{k}}{s_{k}} .
  \end{eqnarray*}
  Then $\ds c_{k}s_{k}=\log\log{k}$.  Also let \\
  
  \begin{eqnarray*}
  g_{s}(x)= \frac{\log{\log{\frac{2}{x}}}+1}{\frac{x}{2}(\log{\frac{2}{x}})^{s+1}(\log{\log{\frac{2}{x}}})^{s+1}}\chi_{[0,\epsilon_{s}]} .
  \end{eqnarray*}
  Where $\epsilon$ is chosen sufficiently small so that $\ds g_{s}$ is monotone decreasing and all expressions involving the logarithms are positive and well defined.   \\
  We have that $\ds (\log{g_{s}(x)}^{p})= (\log(\frac{2}{x}))^{p} C_{s}(x)$  where $C_{s}(x)$ is a bounded function.\\
  So 
  \begin{eqnarray*}
   g_{s}(x)(\log{g_{s}(x)}^{p}) = K_{s}(x)\frac{1}{\frac{x}{2}(\log{\frac{2}{x}})^{s-p+1} (\log{\log{\frac{2}{x}}})^{s}}
   \end{eqnarray*}
    for $K_{s}(x)$ bounded.\\
  The resulting function is integrable provided $s-p+1>1$ and if $s-p+1=1$ we must have $s>1$.\\ 
  Hence $g_{s}(x)$ is in $\ds\mbox{LLog}^{p}\mbox{L}$ for $s>p$ if $s\leq1$ and for $s\geq p$ if $s>1$.\\
  Now let 
  \begin{eqnarray*}
   A_{\lambda}=\frac{1}{\lambda} \int_{0}^{\lambda} g_{s}(x) dx = \frac{C}{\lambda (\log{\frac{2}{\lambda}} \log{\log\frac{2}{\lambda}})^{s}}
   \end{eqnarray*}
  If $\ds\lambda < \frac{1}{k\log{k}}$ we have $\ds A_{\lambda} \geq Cs_{k}$ hence $\ds a_{k} > \frac{1}{k\log{k}}$ and $\sum a_{k} = \infty$ .  \\
  Therefore there exists a perturbed block sequence $\bigcup B_{k} \cup D_{k}$ with $d_{k} = c_{k}l_{k}$ such that the averages of $g_{s}(x)$ fail to converge a.e. along this subsequence.  \\
  Since 
  \begin{eqnarray*}
  c_{k} & = & \frac{\log{\log{k}}}{s_{k}} = \frac{(\log{\log{k})^{s+1}} (\log{k})^{s-1}}{k} ,\\
  \sum_{k=1}^{\infty} \frac{1}{\frac{1}{c_{k}} (\log{\frac{1}{c_{k}}})^{p}} & = & \sum_{k=1}^{\infty} \frac{ \phi(k) (\log{\log{k}})^{s+2}}{k (\log{k})^{p+1-s}}
  \end{eqnarray*}  
  where $\phi(k)$ is a bounded sequence.  This sum converges when $p>s$.  Therefore the averages along this sequence converge for functions in $\mbox{LLog}^{p}\mbox{L}$ with $p>s$.  
\end{ex}
           
\end{document}